\documentclass[12pt]{amsart}
\usepackage{amssymb}

\newtheorem{theorem}{Theorem}

\newtheorem{proposition}{Proposition}
\newtheorem{corollary}{Corollary}
\newtheorem{remark}{Remark}
\newtheorem{definition}{Definition}
\newtheorem{lemma}{Lemma}

\newcommand{\T}{{\mathcal T}}

\newcommand{\F}{{\mathcal F}}
\newcommand{\A}{{\mathcal A}}
\newcommand{\LL}{{\mathcal L}}
\newcommand{\K}{{\mathcal K}}
\newcommand{\Z}{{\mathbb Z}}

\newcommand{\R}{{\mathbb R}}

\newcommand{\PP}{{\mathcal P}}

\begin{document}

\title[Chord diagrams, braid groups and degenerate equality conditions]{Varieties of chord diagrams, braid group cohomology and degeneration of equality conditions}

\author{V.A.~Vassiliev}
\address{Weizmann Institute of Science}

 \email{vavassiliev@gmail.com}

\begin{abstract}
For any  finite-dimensional vector space  $\F$  of  continuous  functions  $f:\R^1 \to \R^1$ consider subspaces in $\F$ defined by systems of {\em equality conditions} $f(a_i) = f(b_i)$, where $(a_i, b_i)$, $i=1, \dots, n$, are some pairs of points in $\R^1$. It is proved that if $\dim \F < 2n-I(n)$, where $I(n)$ is the number of ones in the binary notation of $n$, then there necessarily are independent systems of $n$ equality conditions defining the subspaces of codimension greater than $n$ in $\F$. We also prove lower estimates of the sizes of the inevitable drops of the codimensions of these subspaces.

These estimates are then applied to knot theory (in which systems of equality conditions are known as {\em chord diagrams}). The inevitable presence of complicated non-stable terms  in sequences of spectral sequences calculating the cohomology groups of spaces of knots is proved. 
\end{abstract}



\maketitle

\section{Main results}

Let $\F^N$ be a $N$-dimensional vector subspace of the space $C^0(\R^1,\R^1)$ of continuous functions $\R^1 \to \R^1$. As a rule, a collection of $n $ independent conditions of the form 
\begin{equation}
f(a_i) = f(b_i),
 \label{meq}
\end{equation}
$a_i \neq b_i$,
 $i=1, \dots, n$, defines a subspace of codimension $n$ in $\F^N$ if $n \leq N$, and only the trivial subspace if $n \geq N$. However, for exceptional sets of such conditions the codimensions of these subspaces can drop.

For example, if $\F^N$ is the space $\PP^N$ of all polynomials of the form
\begin{equation} 
\label{ppn}
\alpha_1 x^N + \alpha_2 x^{N-1} + \cdots + \alpha_N x
\end{equation}
in variable $x$, then all subspaces defined by arbitrarily many conditions $f(a_i)=f(-a_i)$ contain the $[N/2]$-dimensional subspace of even polynomials. Of course, the case of polynomials is very specific, but the situation when the dimensions of subspaces in  $\F^N$  defined by some $n$ independent conditions (\ref{meq}) is greater than $\max(N-n,0)$ can be unavoidable by a choice of space $\F^N$.

\begin{definition} \rm
An unordered pair $(a,b)$ of distinct points in $ \R^1$ is called a {\it chord}. An unordered collection of $n$ pairwise distinct chords is called a $n$-{\em chord diagram}. 

The {\em subalgebra of $C^0(\R^1,\R^1)$ corresponding to the $n$-chord diagram}
\begin{equation}
\label{cd}
\{(a_1,b_1), \dots, (a_n,b_n)\}
\end{equation}
 consists of all functions satisfying $n$ conditions (\ref{meq}) with these $a_i, b_i$. These $n$ conditions (and corresponding $n$ chords) are {\em independent} if the codimension of this subalgebra in $C^0(\R^1,\R^1)$ is equal to $n$.
(We say that an affine or vector subspace $\T$ of a function space $\K$ has  codimension $n$ if for any point  $\varphi \in \T$  there exist  $n$-dimensional affine subspaces in  $\K$  intersecting  $\T$  at this point only, and all affine subspaces of higher dimensions through $\varphi$  intersect $\T$ along subspaces of positive dimensions.)
Two independent $n$-chord diagrams are {\em equivalent} if the corresponding subalgebras in $C^0(\R^1,\R^1)$ coincide. A {\em resonance} of a chord diagram
is a cyclic sequence of $k \geq 3$ its pairwise different chords such that one point of each chord belongs also to the preceding chord in this sequence, and the other its point belongs also to the next chord.
\end{definition}

For example, two chord diagrams are equivalent if one of them contains chords $(a, b)$ and $(b, c)$, the other one contains chords $(a, b)$ and $(a, c)$, and all other chords in them are common. 

\begin{proposition}
An $n$-chord diagram is independent if and only if it does not contain resonances.
Two independent $n$-chord diagrams are equivalent if and only if they can be connected by a chain of elementary flips described in the previous paragraph. The space of independent $n$-chord diagrams is a $(2n)$-dimensional smooth manifold.
\end{proposition}

\noindent
{\it Proof} is elementary. \hfill $\diamond$

\subsection{Results for the case of $N \geq n$}

\begin{proposition}
\label{pro}
If $N \geq 2n-1$, then the codimension of the subspace in the space $\PP^N$ of polynomials $($\ref{ppn}$)$, defined by $n$ conditions $($\ref{meq}$)$ of an arbitrary independent $n$-chord diagram $($\ref{cd}$)$, is equal to $n$.
\end{proposition}

\noindent
{\it Proof}. First, the assertion of our proposition will be true if we replace in it the space $\PP^N$ by the $(N+1)$-dimensional space $\widehat{\PP}^N$ of all polynomials of degree $N$. Indeed, any $n$-chord diagram has at most $2n$ distinct endpoints $a_i, b_i$, therefore by interpolation theorem the evaluation morphism from the space of such polynomials to the space of real-valued functions on the set of these endpoints is epimorphic, and hence the preimage of any subspace of codimension $n$ of the latter space also has codimension $n$ in $\widehat{\PP}^N$. However, adding the constant functions preserves the subspace of $\widehat{\PP}^N$ defined by any chord diagram, therefore the codimension of the considered subspace in $\PP^N$ is also equal to $n$. \hfill $\diamond$ \medskip

Denote by $I(n)$ the number of ones in the binary notation of $n$. 

\begin{theorem} \label{mthm}
If $n \leq N< 2n-I(n)$, then for any $N$-dimensional vector subspace $\F^N \subset C^0(\R^1,\R^1)$ there exist independent $n$-chord diagrams $($\ref{cd}$)$, for any of which the codimension of the subspace in $\F^N$ consisting of functions satisfying all the corresponding conditions $($\ref{meq}$)$ is less than $n$. The dimension of the set of such exceptional $n$-chord diagrams is at least $3n-N-1$ in the following exact sense: there exist non-trivial elements of the $(N-n+1)$-dimensional homology group of the $2n$-fold of all independent $n$-chord diagrams, such that each cycle representing such an element necessarily intersects our set.
\end{theorem}

In particular, if $n$ is a power of 2 then the minimal dimension of function spaces $\F^N$ in which any independent $n$-chord diagram defines a subspace of codimension exactly $n$ is equal to $2n-1$. 

 A more general result can be formulated in the terms of configuration spaces; see, for example, \cite{CC} for the current state of the theory of these spaces.

\begin{definition} \rm
The {\it $n$-th configuration space}  $B(X,n)$  of a topological space  $X$  is the (naturally topologized) space of unordered subsets of cardinality  $n$  in  $X$. The {\it regular bundle}  $\xi_n$  with base  $B(X, n)$  is the vector bundle, whose fiber over a  $n$-point configuration is the space of real-valued functions on the corresponding set of points. 
\end{definition}

\begin{theorem}
\label{mthm2}
Suppose that $N \geq n$ and for some natural $r$ the cohomological product 
\begin{equation}
\label{proddd}
 \prod_{i=1}^r w_{N-n+2i-1}(\xi_n)
\end{equation}
 of Stiefel--Whitney classes of the regular bundle $\xi_n$ is not equal to 0 in the ring  $H^*(B(\R^2,n),\Z_2)$. Then for any $N$-dimensional vector subspace  $\F^N$  of the space  $C^0(\R^1,\R^1)$  there exist independent systems of $n$ conditions $($\ref{meq}$)$ such that all subspaces of  $\F^N$  defined by these systems have codimensions  $\leq n-r$ in  $\F^N$.
\end{theorem}

The first statement of Theorem \ref{mthm} follows immediately from this theorem (the case $r=1$) and statement 5.3 of \cite{fuks} asserting that classes $w_{k} \in H^{k}(B(\R^2,n), \Z_2)$ are non-trivial for all $k \leq n-I(n)$; see also Proposition \ref{fu2} in \S \ref{corolpr} below. The second statement of Theorem \ref{mthm} will be proved at the end of \S \ref{pml}.

\begin{corollary} \label{corr}
A. If two natural numbers $n$ and $N$ satisfy one of the following pairs of conditions:

\begin{enumerate}
\item
$n \geq 6, N= n$, 

\item 
$ n\geq 10, N = n+1$, 

\item 
$n\geq 14, N = n+2$ or $n+3$, 

\item
$ n \geq 16, N = n+4 $, 

\item
$ n\geq 18, N = n+5 $, 

\item
$n \geq 20, N = n+6 $,

\item
$n \geq 24, N = n+7$,

\item 
$n \geq 28, N= n+8$ or $n+9$,

\item
$n \geq 32$, $N=n+10$ or $n+11$,
\end{enumerate}

\noindent
then for any $N$-dimensional vector subspace $\F^N \subset C^0(\R^1,\R^1)$ there exist systems of $n$ independent conditions $($\ref{meq}$)$ defining subspaces of codimension $\leq n-2$ in $\F^N$.

B. If $n$ and $N$ satisfy one of the following pairs of conditions:

\begin{enumerate}
\item
$ n \geq 18, N = n $ or $n+1,$ 

\item
$n \geq 22, N = n+2,$ 

\item
$ n\geq 26, N = n+3$, 

\item
$n \geq 30, N = n+4$, 

\item
$ n \geq 36, N=n+5$,

\item
$n \geq 40$, $N=n+6$ or $n+7$,

\item
$n \geq 44$, $N=n+8$ or $n+9$,
\end{enumerate}

\noindent
 then for any $N$-dimensional subspace $\F^N \subset C^0(\R^1,\R^1)$ there exist systems of $n$ independent conditions $($\ref{meq}$)$ defining subspaces of codimension $\leq n-3$ in $\F^N$.

C. If $n$ and $N$ satisfy one of the following conditions: 
\begin{enumerate}
\item $n \geq 30$ and $N = n$ or $n+1$, 

\item $n \geq 44 $ and $N=n+2$ or $n+3$,

\item $n \geq 52 $ and $N = n+4$ or $n+5$,

\item $n \geq 56 $ and $N=n+6$ or $n+7$,
\end{enumerate} 

\noindent
then for any $N$-dimensional subspace $\F^N \subset C^0(\R^1,\R^1)$ there exist systems of $n$ independent conditions $($\ref{meq}$)$ defining subspaces of codimension $\leq n-4$ in $\F^N$.

D. If $n$ and $N$ satisfy one of the following conditions:
\begin{enumerate}
\item 
$n \geq 48$ and $N = n$ or $n+1$, 

\item
$n \geq 60 $ and $N=n+2 $ or $n+3$, 

\item $n \geq 68$ and $N=n+4$ or $n+5$,
\end{enumerate}
then for any $N$-dimensional subspace $\F^N \subset C^0(\R^1,\R^1)$ there exist systems of $n$ independent conditions $($\ref{meq}$)$ defining subspaces of codimension $\leq n-5$ in $\F^N$.

E. If $n$ and $N$ satisfy one of the following conditions:
\begin{enumerate}
\item $n \geq 64$ and $N = n$ or $n+1$, 

\item $n \geq 76 $ and $N=n+2$ or $n+3$,
\end{enumerate}
then for any $N$-dimensional subspace $\F^N \subset C^0(\R^1,\R^1)$ there exist systems of $n$ independent conditions $($\ref{meq}$)$ defining subspaces of codimension $\leq n-6$ in $\F^N$.

F. If $n \geq 80$ and $N = n$ or $n+1$, then for any $N$-dimensional subspace $\F^N \subset C^0(\R^1,\R^1)$ there exist systems of $n$ independent conditions $($\ref{meq}$)$ defining subspaces of codimension $\leq n-7$ in $\F^N$.
\end{corollary}

See \S \ref{corolpr} for the proof of this corollary. Its lists can easily be continued, and the corresponding calculations can be programmed. 

\begin{remark} \rm
The first statement of Theorem \ref{mthm} looks very similar (and is closely related) to the result of \cite{CH} estimating the dimensions of the spaces of functions $\R^2 \to \R^1$ realizing $n$-regular embeddings of the plane. The main effort of our proof of Theorem \ref{mthm2} is a comparison of configuration spaces used in these two problems, see Lemma \ref{le3} below.
\end{remark}

\subsection{Results for the case of $N \leq n$}

\begin{theorem}
\label{mthm3}
If $N \leq n$ and for some natural $r$ the product \begin{equation}
\label{proddd2}
 \prod_{i=1}^r w_{n-N+2i-1}(\xi_n)
\end{equation}
 of Stiefel--Whitney classes of the bundle $\xi_n$ is not equal to 0 in $H^*(B(\R^2,n), \Z_2)$, then for any $N$-dimensional vector subspace $\F^N \subset C^0(\R^1,\R^1)$ 
there exist independent $n$-chord diagrams, for any of which the subspace of $\F^N$ consisting of functions satisfying the corresponding system of equality conditions is at least $r$-dimensional.
\end{theorem}

If $N=n$, then Theorems \ref{mthm2} and \ref{mthm3} coincide tautologically.

\begin{corollary}
\label{cor25}
If $N \geq 2$, then for any $N$-dimensional vector subspace $\F^N \subset C^0(\R^1,\R^1)$ there exist independent $n$-chord diagrams with arbitrarily large $n$ such that the corresponding systems of equality conditions have non-trivial solutions in $\F^N$. 
\end{corollary}

\noindent
Indeed, it is enough to prove this for $N=2$ and numbers $n$ equal to powers of 2. In this case $w_{n-N+1}(\xi_n) \neq 0$ by the previously mentioned result of \cite{fuks}. \hfill $\diamond$

\begin{remark} \rm
This corollary has also an elementary proof. Indeed, any 2-dimensional subspace of $C^0(\R^1,\R^1)$ contains a non-zero function taking equal values at some two different points $a, b \in \R^1$. Then this function necessarily satisfies the equality conditions $f(\tilde a) = f(\tilde b)$ for a continuum of different pairs $(\tilde a, \tilde b) \subset [a, b]$.
\end{remark}

\begin{corollary}
\label{cor26}
All statements of Corollary \ref{corr} will remain valid, if in each of its conditions we replace the value of $N$ by $2n-N$ $($e.g., $N=n+4$ by $N=n-4)$ and simultaneously the corresponding conclusion ``there exist systems of $n$ independent conditions $($\ref{meq}$)$ defining subsets of codimension $\leq n-r$ in $\F^N$'' by ``there exist systems of $n$ independent conditions $($\ref{meq}$)$ defining subsets of dimension $\geq r$ in $\F^N$''.
\end{corollary}

\begin{remark} \rm
In terms of \cite{cs}, the subspaces of abnormal codimensions defined by chord diagrams in finite-dimensional function spaces are responsible for the non-stable regions of the $(p,q)$-planes of spectral sequences converging to cohomology groups of spaces of long knots $\R^1 \to \R^3$ defined by functions from these function spaces. These domains are the only possible source of cohomology classes of the knot space (including 0-dimensional classes, i.e. knot invariants) not of finite type. In \S \ref{alex} below, we prove some facts on filtrations of simplicial resolutions of discriminant spaces in finite-dimensional knot spaces, estimating the deviation of the corresponding spectral sequences from stable ones.
\end{remark}

\section{Scheme of proof of Theorem \ref{mthm2} }
\label{scheme}

Denote by $\mbox{CD}_n$ the set of equivalence classes of independent $n$-chord diagrams. It has a natural topology induced by the topology of the variety of subalgebras of codimension $n$ in $C^0(\R^1,\R^1)$. To describe this topology without infinite-dimensional considerations, let $\A$ be a sufficiently large finite-dimensional vector subspace of $C^0(\R,\R)$, such that all subspaces of $\A$ defined by independent $n$-chord diagrams (that is, the intersections of $\A$ with subalgebras of $C^0(\R^1,\R^1)$ corresponding to these chord diagrams) have codimension exactly $n$ in $\A$, and non-equivalent $n$-chord diagrams define different subspaces. (For reasons similar to the proof of Proposition \ref{pro} we can take for such a space  $\A$  the space  $\PP^N$, $N \geq 2n+1$, or any space containing it; taking $\PP^{2n-1}$ is not enough because non-equivalent $n$-chord diagrams can define equal subspaces in it). We can and will assume that $\A$ contains $\F^N$, because otherwise we can replace  $\A$  by its linear hull with  $\F^N$. 

The set $\mbox{CD}_n$ is embedded into the Grassmann manifold $G(\A,-n)$ of subspaces of codimension $n$ in $\A$, and inherits a topology from this manifold. It is easy to see that this definition of a topology on $\mbox{CD}_n$ does not depend on the choice of $\A$. 

Suppose that $N \geq n$. Let $\Delta_r(\F^N) \subset G(\A,-n)$ be the set of all subspaces of codimension  $n$  in  $\A$ whose linear hulls with  $\F^N$ have codimensions at least  $r$  in $\A$.

\begin{proposition}
The class in $H^*(G(\A,-n),\Z_2)$ Poincar\'e dual to the homology class of algebraic variety $\Delta_r(\F^N)$ 
is equal to $r \times r$ determinant
\begin{equation}
\label{arr}
\begin{array}{|cccccc|}
w_{N-n+r} & w_{N-n+r+1} & w_{N-n+r+2} & \dots & w_{N-n+2r-2} & w_{N-n+2r-1} \\
w_{N-n+r-1} & w_{N-n+r} & w_{N-n+r+1} & \dots & w_{N-n+2r-3} & w_{N-n+2r-2} \\
w_{N-n+r-2} & w_{N-n+r-1} & w_{N-n+r} & \dots & w_{N-n+2r-4} & w_{N-n+2r-3} \\
\dots & \dots & \dots & \dots & \dots & \dots \\
w_{N-n+2} & w_{N-n+3} & w_{N-n+4} & \dots & w_{N-n+r} & w_{N-n+r+1} \\
w_{N-n+1} & w_{N-n+2} & w_{N-n+3} & \dots & w_{N-n+r-1} & w_{N-n+r} 
\end{array}  \ , 
\end{equation}
where $w_i$ are Stiefel--Whitney classes of the tautological bundle over $G(\A,-n)$.
\end{proposition}

\noindent
{\it Proof.} Let $\tau_n$ be the vector bundle over $G(\A,-n)$, whose fiber over the point $\{L\}$ corresponding to subspace $L \subset \A$ is the space of linear forms on $\A$ vanishing on $L$. Consider the morphism of this bundle to the constant bundle with fiber $(\F^N)^*$, sending any such linear form to its restriction to $\F^N$. The variety $\Delta_r(\F^N)$ can be redefined as the set of points $\{L\}$ such that the rank of this morphism does not exceed $n-r$. 
By the real version of Thom--Porteous formula (the proof of which literally repeats its complex analog given in \cite{ful}, \S 14.4, after standard replacements of Chern classes by Stiefel--Whitney classes, $\Z$ by $\Z_2$, etc.) the class in $H^*(G(\A, -n), \Z_2)$ Poincar\'e dual to this variety is equal to the determinant of the form (\ref{arr}) in which all symbols $w_i$ are Stiefel--Whitney classes of the virtual bundle $ - \tau_n$. 

The constant bundle on $G(\A,-n)$ with fiber $\A^*$ is obviously isomorphic to the direct sum of $\tau_n$ and the bundle dual (and hence isomorphic) to the tautological bundle. Therefore $- \tau_n$ and this tautological bundle belong to the same class of the group $\tilde K(G(\A,-n))$, in particular have the same Stiefel--Whitney classes. \hfill $\diamond$
 \medskip

These Stiefel-Whitney classes $w_i(-\tau_n)$ are equal to $i$-dimensional components $\bar w_{i}(\tau_n) \in H^{i}(G(\A,-n), \Z_2)$ of the class $w^{-1}(\tau_n)$, where $w(\tau_n) = 1 + w_1(\tau_n)+ \dots $ is the total Stiefel-Whitney class of bundle $\tau_n$, see \S 4 of \cite{MS}. 
 If the intersection of the subset $\mbox{CD}_n \subset G(\A, -n)$ with $\Delta_r(\F^N)$ is empty then the restriction of the class (\ref{arr}) is trivial in the cohomology group of $\mbox{CD}_n$. Therefore, Theorem \ref{mthm2} is reduced to the following lemma.

\begin{lemma} \label{mlem}
If the class $($\ref{proddd}$)$ is not equal to 0, then 
the restriction of the class $($\ref{arr}$)$ to the subvariety $\mbox{\rm CD}_n \subset G(\A,-n)$ is a non-trivial element of the group $ H^{r(N-n+r)}(\mbox{\rm CD}_n, \Z_2)$.
\end{lemma}

\section{Proof of Lemma \ref{mlem}}
\label{pml}

Let $\R^2_+ \subset \R^2$ be the half-plane $\{(a,b)| a<b\} \subset \R^2$. Any point $(a,b)$ of $\R^2_+$ can be identified with the chord $(a,b)$, and any element of the configuration space $B(\R^2_+,n)$ with a $n$-chord diagram. Denote by $\Xi \subset B(\R^2_+,n)$ the set of dependent (that is, containing resonances) $n$-chord diagrams. 
Consider the diagram 
\begin{equation} \label{diag}
\begin{array}{ccccc}
B(\R^2,n)& \supset & B(\R^2_+,n) \setminus \Xi & & \\
& &  \ \pi {\downarrow} & & \\
 & &  \  \ \mbox{CD}_n & \subset &G(\A,-n)  ,
\end{array}
\end{equation} 
where $\pi$ is the map sending any chord diagram to its equivalence class.

\begin{lemma}\label{le1}
The restriction of the regular vector bundle $\xi_n$ to the  subset $B(\R^2_+,n) \setminus \Xi \subset B(\R^2,n)$ is isomorphic to the bundle pulled back by the map $\pi$ from the bundle $\tau_n$ over $\mbox{CD}_n$. 
\end{lemma}

\noindent 
{\it Proof.} The bundle $\tau_n$ is  isomorphic to its dual bundle $\tau_n^*$, i.e. the quotient of the trivial bundle with fiber $\A$ by the tautological bundle over $G(\A, -n)$.

Consider the following homomorphism from the trivial bundle with fiber $\A$ over $B(\R^2_+,n) \setminus \Xi $ to $\xi_n$: over any $n$-chord diagram $\Gamma$ it sends any function $f \in \A \subset C^0(\R^1,\R^1)$ to the function on the set of chords of this chord diagram, whose value on any chord $(a_i,b_i)$ is equal to the difference $f(b_i)-f(a_i)$. By the first characteristic property of the space $\A$ this morphism is surjective; by definition of inclusion 
$\mbox{CD}_n \subset G(\A,-n)$
 its kernel is equal to the fiber of the tautological bundle over the point $\pi(\Gamma) \in \mbox{CD}_n $. Therefore our homomorphism induces an isomorphism between the bundles $\pi^*(\tau^*) \sim \pi^*(\tau)$ and $\xi_n$. \hfill $\diamond$ 

\begin{lemma}[see \cite{fuks} or Proposition \ref{fu3} below] \label{le2}
 The square of any positive-dimensional element of the ring $H^*(B(\R^2,n), \Z_2)$ is equal to zero, in particular 
$w^{-1}(\xi_n) = w(\xi_n)$ and 
$w_i(\xi_n) = \bar w_i(\xi_n)$ for any $i$. \hfill $\diamond$
\end{lemma}

\begin{lemma}
\label{lem9}
The determinant of the form $($\ref{arr}$)$ in which all classes $w_i$ are replaced by
 $w_i(\xi_n)$ is equal to the product $($\ref{proddd}$)$ in $H^*(B(\R^2,n),\Z_2)$. 
\end{lemma}

\noindent
{\it Proof.} The matrix (\ref{arr}) is symmetric with respect to the side diagonal, hence calculating its determinant $\mbox{mod} \, 2$ it is enough to count only those products of $r$ matrix elements which are self-symmetric with respect to this diagonal. By Lemma \ref{le2} such products, not all factors of which lie in this diagonal, are also trivial. \hfill $\diamond$

\begin{lemma}\label{le3}
The inclusion $B(\R^2_+,n) \setminus \Xi \to B(\R^2,n)$ induces a monomorphism of cohomology groups $H^*(B(\R^2,n),\Z_2) \to H^*(B(\R^2_+,n) \setminus \Xi, \Z_2).$
\end{lemma} 

\noindent
Lemma \ref{le3} will be proved in \S \ref{pl4}. Lemma \ref{mlem} follows from Lemmas \ref{le1}--\ref{le3} and the functoriality of Stiefel--Whitney classes. Namely, by Lemma \ref{le3} if the product (\ref{proddd}) is non-trivial in $H^*(B(\R^2,n), \Z_2)$, then it is non-trivial also in $H^*(B(\R^2_+,n) \setminus \Xi, \Z_2)$. By Lemmas \ref{le1}--\ref{le3} this element of $H^*(B(\R^2_+,n) \setminus \Xi, \Z_2)$ equals the class induced by the map $\pi$ from determinant (\ref{arr}), hence this determinant is also non-trivial.
\medskip

\noindent{\it Proof of the last statement of Theorem \ref{mthm}.} By Lemma \ref{le3} and statement 5.3 of \cite{fuks}, the class $w_{N-n+1}(\xi_n)$ is non-trivial. Then we can take an arbitrary element of the group $H_{N-n+1}(B(\R_+^2,n) \setminus \Xi, \Z_2)$ on which this class takes non-zero value: any cycle realizing such an element intersects the set $\pi^{-1}(\Delta_1(\F^N))$. \hfill $\diamond$

\section{Proof of Theorem \ref{mthm3}}

Now suppose that $N \leq n$. Let $\Lambda_r(\F^N)$ be the subset of $G(\A,-n)$ consisting of planes, whose intersection with $\F^N$ is at least $r$-dimensional.

\begin{proposition}
\label{propz}
The class in $H^*(G(\A,-n), \Z_2)$ Poincar\'e dual to the variety $\Lambda_r(\F^N)$ 
is equal to $r \times r$ determinant similar to $($\ref{arr}$)$, in which $N-n$ in all lower indices is replaced by $n-N$, and $w_i$ are Stiefel--Whitney classes of the bundle $\tau_n$.
\end{proposition}

\noindent
{\it Proof.}
The projection along the fibers of the tautological bundle over $G(\A,-n)$ defines a
morphism from the constant bundle with fiber  $\F^N$ and base $G(\A, -n)$ to the bundle dual (and hence isomorphic) to $\tau_n$, i.e. to the quotient of the constant bundle with fiber $\A$ by the tautological bundle. The set $\Lambda_r(\F^N)$ can be defined as the set of points, at which the rank of this morphism does not exceed $N-r$. Our proposition follows again from the real version of Thom--Porteous formula applied to this morphism. \hfill $\diamond$ \medskip

The rest of the reduction of Theorem \ref{mthm3} to Lemma \ref{le3} repeats that of Theorem \ref{mthm2}; the difference between Stiefel--Whitney classes of the bundles $-\tau_n$ and $\tau_n$ participating in the corresponding Thom--Porteous formulas is eliminated by Lemma \ref{le2}.

\section{Proof of Lemma \ref{le3}} 
\label{pl4}

We will prove the dual statement: the map $H_*(B(\R^2_+,n) \setminus \Xi, \Z_2) \to H_*(B(\R^2,n),\Z_2)$ induced by identical embedding is epimorphic. 

By \cite{fuks}, all stabilization maps $H_*(B(\R^2,n),\Z_2) \to H_*(B(\R^2,n+m),\Z_2) $ induced by standard inclusions $B(\R^2,n) \hookrightarrow B(\R^2,n+m)$ are injective. Therefore, all elements of the group $H_*(B(\R^2,n),\Z_2)$ are given by polynomials in multiplicative generators of Hopf algebra $H_*(B(\R^2,\infty), \Z_2)$, and it is enough to prove that all these generators and their products participating in the construction of these elements can be realized by cycles lying in $B(\R^2_+, n) \setminus \Xi$. 

These generators $[M_j] \in H_{2^j-1}(B(\R^2,2^j),\Z_2)$ were defined in \S 8 of \cite{fuks} by the following cycles $M_j \subset B(\R^2,2^j)$. Choose somehow two opposite points of the circle of radius $1$ with center at the origin in $\R^2$. Take two circles of a small radius $\varepsilon$ with centers at these points and choose somehow two opposite points in any of them. Take circles of radius $\varepsilon^2$ with centers at all obtained four points and choose two opposite points in any of them... Continuing, we obtain after the $j$-th step a $2^j$-configuration in $\R^2$. This construction includes $1+2+4 + \dots +2^{j-1}$ choices of opposite points in certain circles, so the set $M_j$ of all possible $2^j$-configurations which can be obtained in this way is $(2^j-1)$-dimensional. It is easy to see that this set is a closed submanifold in $B(\R^2,2^j)$, and so it defines an element of the group $H_{2^j-1}(B(\R^2,2^j),\Z_2)$, $j \geq 1$. Finally, define the element $[M_0]\in H_0(B(\R^2,1), \Z_2)$ as the class of a single point.

Unfortunately, these cycles with $j>1$ contain configurations with resonances, and, moreover, all configurations of $M_j$ do not lie in $\R^2_+$. To avoid these problems, we modify the previous construction by 1) replacing circles with squares, 2) taking these squares of varying sizes, and 3) shifting the resulting configurations into the half-plane $\R^2_+ \subset \R^2$; see the next subsection. 

\subsection{Construction of cycles $\tilde M_j$ (see Fig. \ref{conf})}

 Let us fix a very small number $\varepsilon >0$. Define the {\it basic square} $\Box \subset \R^2$ as the union of four segments
connecting consecutively the points $(-1,-1), (-1,1), (1,1), (1,-1)$ and again $(-1,-1)$. Fix an arbitrary continuous function $\chi: \Box \to [\varepsilon,1]$ equal identically to 1 on segments 
\begin{equation}\label{mac}
 [(-1+\varepsilon,1),(1,1)] \qquad \mbox{and} \qquad [(1,-1+\varepsilon),(1,1)],
\end{equation}
equal to $\varepsilon$ on segments 
\begin{equation}\label{mic}
 [(-1,-1),(-1,1-\varepsilon)] \qquad \mbox{and} \qquad [(-1,-1),(1-\varepsilon,-1],
\end{equation}
and taking some intermediate values in remaining $\varepsilon$-neighborhoods of corners $(-1,1)$ and $(1,-1)$. 

Further, define two sequences of natural numbers $u_j$ and $T_j$, $j \geq 2$, by the recursion 
\begin{equation} \label{fib}
u_2=1, T_2=2; \quad u_j=u_{j-1}+T_{j-1}+2, \quad T_j = u_j+ T_{j-1} +1 \quad \mbox{for}  j>2.
\end{equation}

\unitlength=0.7mm

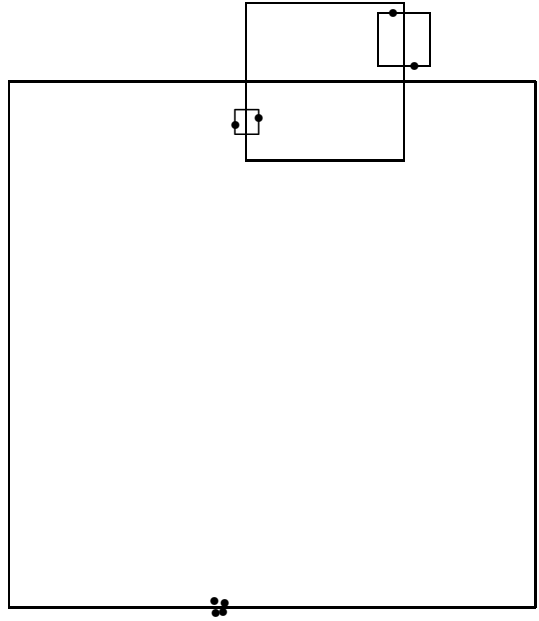
\begin{figure}
\begin{picture}(120,117)
\put(10,2){\line(1,0){100}}
\put(10,2){\line(0,1){100}}
\put(110,102){\line(-1,0){100}}
\put(110,102){\line(0,-1){100}}
%\put(70,102){\circle{1.5}}
%\put(50,2){\circle{1.5}}
\put(55,87){\line(1,0){30}}
\put(55,87){\line(0,1){30}}
\put(85,117){\line(-1,0){30}}
\put(85,117){\line(0,-1){30}}
%\put(55,94){\circle{1.5}}
%\put(85,110){\circle{1.5}}

\put(80,105){\line(1,0){10}}
\put(80,105){\line(0,1){10}}
\put(90,115){\line(-1,0){10}}
\put(90,115){\line(0,-1){10}}

\put(87,105){\circle*{1.5}}
\put(83,115){\circle*{1.5}}
\put(53,93.7){\circle*{1.5}}
\put(57.5,95){\circle*{1.5}}

\put(51,2.9){\circle*{1.5}}
\put(50.7,1.2){\circle*{1.5}}
\put(49.3,1){\circle*{1.5}}
\put(49,3.2){\circle*{1.5}}

\put(52.4,91.8){\large $\Box$}
\end{picture}
\caption{An 8-configuration of class $\tilde M_3$}
\label{conf}
\end{figure}

Define the subset $\tilde M_1 \subset B(\R^2,2)$ as the space of all choices of two opposite points in the basic square $\Box$ .
Suppose that we have defined the variety $\tilde M_{j-1} \subset B(\R^2,2^{j-1})$, $j \geq 2$. The variety $\tilde M_j \subset B(\R^2,2^j)$ is then defined as the space of all $2^j$-configurations consisting of 

1) a $2^{j-1}$-configuration obtained from some configuration of class $\tilde M_{j-1}$ by the affine map $\R^2 \to \R^2$ given by the formula
\begin{equation}\{X \mapsto A + \varepsilon \chi^{u_j}(A)X\} , \label{nine} \end{equation} where $A$ is some point of the basic square $\Box$, and

2) a $2^{j-1}$-configuration obtained from some configuration of class $\tilde M_{j-1}$ by the map 
\begin{equation}\{X \mapsto -A+ \varepsilon \chi^{u_j}(-A) X\}\label{ten}\end{equation}
 with the same $A$.

Any $2^j$-configuration of the class $\tilde M_j$ determines uniquely the set of $2^j-1$ {\it squares participating in its construction}: it consists of the basic square $\Box$ and the images under maps (\ref{nine}), (\ref{ten}) 
 of two collections of $2^{j-1}-1$ squares participating in the construction of two $ 2^{j-1}$-configurations of class $\tilde M_{j-1}$. This set is obviously organized in a binary tree, whose leaves are related with some pairs of points of our configuration.

\begin{lemma} \label{est1}
The length of the sides of the smallest square participating in the construction of a $2^j$-configuration of class $\tilde M_j$ is not less than $2\varepsilon^{T_j}$.
\end{lemma}

\noindent
{\it Proof.} 
This 
follows by induction from the last condition (\ref{fib}) because both subsets constituting our $2^j$-configuration are obtained from certain $2^{j-1}$-configurations of class $\tilde M_{j-1}$ by homotheties (\ref{nine}), (\ref{ten}) with coefficients $\geq \varepsilon^{u_j+1}$. \hfill $\diamond$ \medskip

Finally, we move the obtained variety $\tilde M_j$ into $B(\R^2_+,{2^j})$ shifting all its $2^j$-con\-fi\-gu\-rations to $\R^2_+$ by some translation $\{X \mapsto X+Z\}$, where $Z$ is an arbitrary vector $(u,v) \in \R^2$ with $v-u \geq 8$. Denote by $\nabla_j$ the resulting cycle in $B(\R^2_+,2^j)$.

It is easy to see that $\nabla_j$ is the image of an embedding $M_j \to B(\R^2,2^j)$ which is homotopic to the identical embedding; in particular, it defines the same homology class in $H_*(B(\R^2,2^j), \Z_2)$. Indeed, such a homotopy is defined by 1) a homotopy connecting the function $\chi$ with the function equal identically to 1 in the space of positive functions $\Box \to \R^1$, 2) deformations of all circles to squares, and 3) the continuous shift of the plane by the vector $Z$. Therefore, it remains to prove the following statement.

\begin{theorem}
\label{mtool}
The variety $\nabla_j$ does not contain $2^j$-chord diagrams with resonances.
\end{theorem}

\noindent
{\it Proof.} All chords $(a,b) \in \R^2_+,$ $a<b$, which can occur in some configuration from $\nabla_j$, belong to the $\sqrt{2}/(1-\varepsilon)$-neighborhood of the point $Z \in \R^2_+$, therefore the left-hand point $a_k$ of either of these $n$ chords cannot coincide with the right-hand point $b_l$ of some other chord. Hence, if some configuration $\Gamma \in \nabla_j$ contains a resonance, then any closed chain of its chords defining this resonance consists of an even number of chords $(a_i, b_i) \in \R^2_+$, such that the closed sequence of segments connecting these chords consists of strictly alternating horizontal and vertical segments.

\begin{definition} \rm
For any positive $\delta$, a segment in $\R^2$ is called {\it $\delta$-horizontal} (respectively, {\it $\delta$-vertical}) if the tangent of the angle between this segment and a horizontal (respectively, vertical) line belongs to the interval $[0,\delta)$.
A $n$-configuration in $\R^2$ is $\delta$-{\it resonant} if there exists a closed chain of strictly alternating $\delta$-vertical and $\delta$-horizontal segments in $\R^2$ such that the endpoints of any of its segments belong to our configuration.
\end{definition}

Theorem \ref{mtool} follows now from the following statement.

\begin{theorem}
\label{mtool2}
 If the number $\varepsilon$ participating in the construction of cycle $\tilde M_j$ is sufficiently small,  then  $\tilde M_j$ does not contain $\varepsilon^{u_j+1}$-resonant $2^j$-configurations.
\end{theorem}

\noindent
{\it Proof.} Suppose that this theorem is proved for all cycles $\tilde M_i$, $i<j$. 
By construction, any $2^j$-configuration $\Gamma \in \tilde M_j$ splits into two subsets of cardinality $2^{j-1}$ with mass centers at some opposite points $A$ and $-A$ of the basic square $\Box$, any of these subsets lying in the $\sqrt{2}\varepsilon/(1-\varepsilon)$-neighborhood of the corresponding point $A$ or $-A$. 

Suppose that our configuration $\Gamma \in \tilde M_j$ is $\varepsilon^{u_j+1}$-resonant. If the entire chain of its points participating in this resonance lies in only one of these two subsets of $\Gamma$, then we get a contradiction with the induction hypothesis for $i=j-1$, because this subset is homothetic to a configuration of the class $\tilde M_{j-1}$, and $\varepsilon^{u_j+1}<\varepsilon^{u_{j-1}+1}$.

So, our chain should contain $\varepsilon^{{u_j}+1}$-vertical or $\varepsilon^{u_j+1}$-horizontal segments, which connect some points from these two subsets. Therefore, the corresponding points $A$ and $-A$ are very close either to central points of opposite horizontal sides of the basic square $\Box$, or to central points of its vertical sides. These two situations can be reduced to each other by the reflection in diagonal $\{x=y\}$ of $\R^2$, so it is enough to consider only the first of them.

Consider a $\varepsilon^{{u_j}+1}$-vertical segment of our resonance chain which has endpoints in both these subsets; let $A_0$ be its endpoint in the upper subset. Starting from $A_0$, our chain travels somehow inside this upper subset and finally leaves it along some other $\varepsilon^{u_j+1}$-vertical segment; let $B_0$ be the upper point of the latter segment.

\begin{lemma} \label{est}
1. The distance between the projections of points $A_0$ and $B_0$ to any horizontal line in $\R^2$ is estimated from above by $7\varepsilon^{u_j+1}$.

2. The distance between projections of $A_0$ and $B_0$ to any vertical line is estimated from below by $\varepsilon^{T_{j-1}+1}$.
\end{lemma}

\noindent
{\it Proof.} 
1. This distance is estimated from above by the sum of a) the maximal difference of the $x$-coordinates of points of the lower subset of our $2^j$-configuration $\Gamma$, b) the difference of the $x$-coordinates of $A_0$ and the other endpoint of the segment of our chain connecting $A_0$ with the lower subset of this configuration, and c) the similar difference for $B_0$. By the construction of $\tilde M_j$, the first of these differences is equal to $2\varepsilon^{{u_j}+1}(1 +O(\varepsilon))$, and the other two are estimated from above by the lengths of these two segments 
(equal to $2(1+O(\varepsilon))$) multiplied by their tangents with the vertical direction, which are estimated by $\varepsilon^{u_j+1}$ since these segments are $\varepsilon^{u_j+1}$-vertical. So, the entire distance is estimated from above by $\varepsilon^{u_j+1}(6+O(\varepsilon)) < 7 \varepsilon^{u_j+1}$.

2. Consider the sequences of decreasing squares participating in the construction of our configuration $\Gamma \in \tilde M_j$, which converge to the points $A_0$ and $B_0$. Let $\Box_k$ be the last common square of these two sequences. By Lemma \ref{est1}, the length of its sides is at least $2 \varepsilon^{T_{j-1}+1}$: indeed, this square is obtained from a square participating in the construction of a $2^j$-configuration of class $\tilde M_{j-1}$ by a homothety with coefficient $\varepsilon \chi(A)$ for some point $A$ from the central part of the upper side of the basic square (where $\chi \equiv 1$). The next two squares in these sequences (or their final points $A_0$ and $B_0$ if $\Box_k$ is a square of the last $j$-th level) are different, therefore these next squares (or points) are centered at (or coincide with) some points of opposite sides of $\Box_k$. These cannot be vertical sides: indeed, in this case the $x$-coordinates of our points $A_0$ and $B_0$ would differ by $2\varepsilon^{T_{j-1}+1} (1+O(\varepsilon))$, which contradicts statement 1 of our lemma, 
because by (\ref{fib}) $\varepsilon^{T_{j-1}+1} \gg \varepsilon^{u_j+1}$. Therefore, these are horizontal opposite sides, hence the difference of their $y$-coordinates is estimated from below by the number $2\varepsilon^{T_{j-1}+1} (1+O(\varepsilon)) > \varepsilon^{T_{j-1}+1}$. The last estimate holds also for the difference of $y$-coordinates of satellite points $A_0$ and $B_0$ of these sides. \hfill $\diamond$

\begin{corollary}
Segment $[A_0, B_0]$ is $\varepsilon^{u_{j-1}+1}$-vertical.
\end{corollary}
 
\noindent
{\it Proof.} By the previous lemma, the absolute value of the tangent of the angle between this segment and the vertical direction is estimated from above by $7\varepsilon^{u_j-T_{j-1}}$, which by (\ref{fib}) is less than $\varepsilon^{u_{j-1}+1}$ (since we can assume that $\varepsilon < 1/7$). \hfill $\diamond$ \medskip

In particular, this segment $[A_0, B_0]$ is not $\varepsilon^{u_j+1}$-horizontal, so $A_0$ and $B_0$ cannot be neighboring points in our $\varepsilon^{u_j+1}$-resonance chain.
Now consider the closed chain of segments in the upper subset of our $\tilde M_j$-configuration $\Gamma$, which consists of the segment $[A_0, B_0]$ and the part of our initial $\varepsilon^{u_j+1}$-resonance chain connecting these two points inside this upper subset of $\Gamma$. This closed chain is a $\varepsilon^{u_{j-1}+1}$-resonance, which contradicts the induction hypothesis over $j$. This contradiction finishes the proof of Theorem \ref{mtool2}, and hence also of Theorem \ref{mtool}. \hfill $\diamond$ \medskip

Finally, any product $ \left[ M_{j_1} \right] \cdot \left[M_{j_2} \right] \cdots \left[M_{j_q} \right]$ of multiplicative generators of the Hopf algebra $H_*(B(\R^2,\infty), \Z_2)$ such that $2^{j_1} + 2^{j_2} + \dots + 2^{j_q} = n$ can be realized by the set of $n$-configurations, some $2^{j_1}$ points of which form a configuration of type $\tilde M_{j_1}$ shifted to $\R^2_+$ along the vector $(0,8)$, some other $2^{j_2}$ points form a configuration of type $\tilde M_{j_2}$ shifted along the vector $(8, 16)$ ... and the last $2^{j_q}$ points form a configuration of type $\tilde M_{j_q}$ shifted along the vector $(8(q-1), 8q)$. There are no vertical or horizontal segments in $\R^2$ connecting points of different such groups of configurations, thus all such $n$-configurations do not contain resonances. This finishes the proof of Lemma \ref{le3} and hence also of Theorems \ref{mthm2} and \ref{mthm3}. \hfill $\diamond$

\section{Proof of Corollary \ref{corr}}
\label{corolpr}

\begin{proposition} \label{monot}
All statements of Corollary \ref{corr} are monotone over $N$: if for some triple of numbers $(n, N, r)$ it is true that for any $N$-dimensional subspace $\F^N \subset C^0(\R^1, \R^1)$ there exist systems of $n$ independent equality conditions defining subspaces of codimension $\leq n-r$ in  $\F^N$, then the same is true for the triple $(n, N-1, r)$. 
\end{proposition}

\noindent
{\it Proof.} Apply the hypothesis of this proposition to an arbitrary $N$-dimensional space containing  $\F^{N-1}$. \hfill $\diamond$ \medskip

Let us recall several results of \cite{fuks} on mod 2 cohomology of spaces $B(\R^2,n)$.

\begin{proposition}[see \cite{fuks}, \S 4.8] 
\label{fu1}
For any $k$, the group $H^k(B(\R^2,n),\Z_2)$ has a canonical basis whose elements are in a one-to-one correspondence with unordered decompositions of the number $n$ into $n-k$ powers of 2. In particular, this group is nontrivial if and only if  $k \leq n - I(n)$.
\end{proposition}

The standard notation for such a basis element is $\langle 2^{l_1}, 2^{l_2}, \dots , 2^{l_{t}} \rangle,$ where $l_1 \geq l_2 \geq \dots \geq l_{t} \geq 1$, $t \leq n-k$: it is the list (in non-increasing order) of all summands of such a decomposition {\em which are strictly greater than 1}. 

Namely, such a basis element of $H^*(B(\R^2,n), \Z_2)$ is defined by the intersection index with the closure of the subvariety in $B(\R^2, n)$ consisting of all $n$-configurations such that there exist $t$ distinct vertical lines in $\R^2$, one of which contains $2^{l_1}$ points of our configuration, some other one contains $2^{l_2}$ of them, etc.

We will also use the shortened notation $\langle 2^{s_1}_{v_1}, 2^{s_2}_{v_2}, \dots, 2^{s_q}_{v_q} \rangle$ for these basis elements, where $s_1 > s_2 > \dots > s_q \geq 1$ and $2^{s_i}_{v_i}$ means $2^{s_i}$ repeated $v_i$ times; if some $v_i$ is here equal to 1 then we write simply $2^{s_i}$ instead of $2^{s_i}_1$.

\begin{proposition}[see \cite{fuks}, \S 5.2] 
\label{fu2}
For any $k<n$ the class $w_k(\xi_n) \in H^{k}(B(\R^2,n), \Z_2)$ is equal to the sum of all basic elements of this group described in the previous proposition. In particular, all classes $w_k(\xi_n)$ with $k \leq n-I(n)$ are not equal to 0.
\end{proposition}

So we have 
\begin{equation}
\label{sw1}
w_1 = \langle 2 \rangle
\end{equation}
\begin{equation}
\label{sw2}
w_2 = \langle 2_2 \rangle
\end{equation}
\begin{equation}
\label{sw3}
w_3 = \langle 2_3 \rangle + \langle 4 \rangle
\end{equation}
\begin{equation}
\label{sw4}
w_4 = \langle 2_4 \rangle + \langle 4, 2 \rangle
\end{equation}
\begin{equation}
\label{sw5}
w_5 = \langle 2_5 \rangle + \langle 4, 2_2 \rangle
\end{equation}
\begin{equation}
\label{sw6}
w_6 = \langle 2_6 \rangle + \langle 4, 2_3 \rangle + \langle 4_2 \rangle
\end{equation}
\begin{equation}
\label{sw7}
w_7 = \langle 2_7 \rangle + \langle 4, 2_4 \rangle + \langle 4_2, 2 \rangle + \langle 8 \rangle
\end{equation}
\begin{equation}
\label{sw8}
w_8 = \langle 2_8 \rangle + \langle 4, 2_5 \rangle + \langle 4_2, 2_2 \rangle + \langle 8, 2 \rangle
\end{equation}
\begin{equation}
\label{sw9}
w_9 = \langle 2_9 \rangle + \langle 4, 2_6 \rangle + \langle 4_2, 2_3 \rangle + \langle 4_3 \rangle + \langle 8, 2_2 \rangle
\end{equation}
\begin{equation}
\label{sw10}
w_{10} = \langle 2_{10} \rangle + \langle 4, 2_7 \rangle + \langle 4_2, 2_4 \rangle + \langle 4_3, 2\rangle + \langle 8, 2_3 \rangle + \langle 8, 4 \rangle
\end{equation}
\begin{equation}
\label{sw11}
w_{11}= \langle 2_{11} \rangle + \langle 4, 2_8 \rangle + \langle 4_2, 2_5 \rangle + \langle 4_3, 2_2 \rangle + \langle 8, 2_4 \rangle + \langle 8, 4, 2 \rangle
\end{equation} 
\begin{equation}
\label{sw12}
w_{12} = \langle 2_{12} \rangle + \langle 4, 2_9 \rangle + \langle 4_2, 2_6 \rangle + \langle 4_3, 2_3\rangle + \langle 4_4 \rangle + \langle 8, 2_5 \rangle + \langle 8, 4, 2_2 \rangle 
\end{equation}

\begin{proposition}[see \cite{fuks}, \S\S 9 and 6] 
\label{fu3}
The cohomological product of two basis elements of the group $H^*(B(\R^2,n), \Z_2)$ having the form $$\langle 2^m, \dots, 2^m, 2^{m-1}, \dots, 2^{m-1}, \dots, 2, \dots, 2 \rangle,$$
where any number $2^i$, $i \in \{1, 2, \dots, m\}$, occurs $p_i$ times in the first factor and $q_i$ times in the second $($and some of numbers $p_i$, $q_i$ can be equal to 0$)$, is equal to
\begin{equation} \label{mul} \prod_{i=1}^m\binom{p_i+q_i}{p_i} \langle 2^m, \dots, 2^m, 2^{m-1}, \dots, 2^{m-1}, \dots, 2, \dots, 2 \rangle  ,
\end{equation}
where any symbol $2^i$ in the angle brackets occurs $p_i+q_i$ times, all binomial coefficients are counted modulo 2, and the entire expression $($\ref{mul}$)$ is assumed to be zero if $(p_m+q_m) 2^m + (p_{m-1}+q_{m-1})2^{m-1} + \cdots + (p_1+q_1) 2 > n$.
\end{proposition}

Now all statements of Corollary \ref{corr} follow immediately from Theorem \ref{mthm2} and the following calculations. \medskip

{\bf A(1).} By (\ref{sw1}), (\ref{sw3}) and (\ref{mul}), $w_1 w_3 = \langle 4, 2 \rangle $, which is non-trivial for $n \geq 6$.

{\bf A(2).} By (\ref{sw2}), (\ref{sw4}) and (\ref{mul}), 
\begin{equation}
\label{a24}
w_2 w_4 = \langle 4, 2_3 \rangle + \langle 2_6 \rangle ,
\end{equation}
which is non-trivial if $n \geq 10$.

{\bf A(3).} By (\ref{sw3}), (\ref{sw5}) and (\ref{mul}), $w_3 w_5 = \langle 4, 2_5 \rangle, $ which is non-trivial if $n \geq 14.$ By (\ref{sw4}), (\ref{sw6}) and (\ref{mul}), 
\begin{equation} \label{a46}
w_4 w_6 = \langle 4_2, 2_4 \rangle + \langle 4_3, 2 \rangle , \end{equation}
 which is also non-trivial if $n \geq 14$.

{\bf A(4).} By (\ref{sw5}), (\ref{sw7}) and (\ref{mul}), 
\begin{equation} 
\label{a57}
w_5 w_7 = \langle 4_3, 2_3 \rangle + \langle 8, 2_5 \rangle + \langle 8, 4, 2_2 \rangle ,
\end{equation} 
which is non-trivial for $n \geq 16$.

{\bf A(5).} By (\ref{sw6}), (\ref{sw8}) and (\ref{mul}), 
\begin{equation} \label{a68}
w_6 w_8 = \langle 2_{14} \rangle + \langle 4, 2_{11} \rangle + \langle 4_2, 2_8 \rangle + \langle 4_3, 2_5 \rangle + \langle 8, 2_7 \rangle + \langle 8, 4_2, 2 \rangle, 
\end{equation} 
which is non-trivial for $n \geq 18$. 

{\bf A(6).} By (\ref{sw7}), (\ref{sw9}) and (\ref{mul}), 
\begin{equation} \label{a79}
w_7 w_9 = \langle 4, 2_{13} \rangle + \langle 4_3, 2_7 \rangle + \langle 8, 2_9 \rangle + \langle 8, 4_3 \rangle,
 \end{equation} 
which is non-trivial for $n \geq 20$.

{\bf A(7).} By (\ref{sw8}), (\ref{sw10}) and (\ref{mul}), 
\begin{equation} \label{a810}
 w_8 w_{10} = \langle 4_2, 2_{12} \rangle + \langle 4_3, 2_9 \rangle + \langle 8, 4, 2_8 \rangle + \langle 8, 4_2, 2_5 \rangle + \langle 8, 4_3, 2_2 \rangle, \end{equation}
which is non-trivial if $ n \geq 24$.

{\bf A(8)}. By (\ref{sw9}), (\ref{sw11}) and (\ref{mul}), 
\begin{equation} \label{a911}
w_9 w_{11} = \langle 
4_3, 2_{11} \rangle + \langle 8, 4, 2_{10} \rangle + \langle 8, 4_3, 2_4 \rangle, 
\end{equation}
which is non-trivial if $ n \geq 28$.

By (\ref{sw10}), (\ref{sw12}) and (\ref{mul}), 
\begin{equation} \label{a1012} w_{10} w_{12} = \langle 4_4, 2_{10} \rangle + \langle 4_5, 2_7 \rangle + \langle 4_6, 2_4 \rangle + \langle 4_7, 2 \rangle + \langle 8, 4, 2_{12} \rangle + \langle 8, 4_4, 2_3 \rangle + \langle 8, 4_5 \rangle , \end{equation} which also is non-trivial if $n \geq 28$.

{\bf A(9)}. The class $w_{14}$ contains summand $\langle 8_2\rangle$, therefore by (\ref{sw12}) and (\ref{mul}) the product $w_{12} w_{14}$ contains summands $\langle 8_2, 4_4 \rangle$ and $\langle 8_3, 4, 2_2 \rangle$, each of which is non-trivial if $n \geq 32$.

{\bf B(1).} By (\ref{sw2}), (\ref{a46}) and (\ref{mul}), 
\begin{equation}
\label{b246}
w_2 w_4 w_6 = \langle 4_2, 2_6 \rangle + \langle 4_3, 2_3 \rangle ,
\end{equation}
 which is non-trivial for $n\geq 18$. By Theorem \ref{mthm2} this calculation proves our statement for $N=n+1$, and the case $N=n$ follows by monotonicity, see Proposition \ref{monot}.

{\bf B(2).} By (\ref{sw3}), (\ref{a57}) and (\ref{mul}), $w_3 w_5 w_7 = \langle 8, 4, 2_5 \rangle,$ which is non-trivial for $n \geq 22$.

{\bf B(3).} By (\ref{a46}), (\ref{sw8}) and (\ref{mul}), 
\begin{equation} \label{b3p} w_4 w_6 w_8 = \langle 4_2, 2_{12} \rangle + \langle 4_3, 2_9 \rangle + \langle 8, 4_2, 2_5 \rangle, \end{equation} which is non-trivial for $n \geq 26$.

{\bf B(4).} By (\ref{sw5}), (\ref{a79}) and (\ref{mul}), 
 \begin{equation} \label{b579}
 w_5 w_7 w_9 = \langle 8, 4, 2_{11} \rangle + \langle 8, 4_3, 2_5 \rangle,
 \end{equation} 
which is non-trivial if $n \geq 30$.

{\bf B(5).} By (\ref{sw6}), (\ref{a810}) and (\ref{mul}), 
\begin{equation} \label{b5p} w_6 w_8 w_{10} = \langle 8, 4, 2_{14} \rangle + \langle 8, 4_3, 2_8 \rangle , \end{equation}
which is non-trivial for $n \geq 36$.

{\bf B(6).} By (\ref{a810}), (\ref{sw12}) and (\ref{mul}),
\begin{equation} \label{b81012} w_8 w_{10} w_{12} = \langle 4_6, 2_{12} \rangle + \langle 4_7, 2_9 \rangle + \langle 8, 4_3, 2_{14} \rangle + \langle 8, 4_5, 2_8 \rangle + \langle 8, 4_6, 2_5 \rangle + \langle 8, 4_7, 2_2 \rangle, 
\end{equation} which is non-trivial for $n \geq 40$. By Theorem \ref{mthm2}, this implies statement B(6) for $N=n+7$, and the case $N=n+6$ follows by monotonicity. Notice that the routine consideration for $N=n+6$ based on formula 
\begin{equation} w_7 w_9 w_{11} = \langle 4_3, 2_{11} \rangle + \langle 8, 4_3, 2_{11} \rangle , \end{equation} gives the same result in more restrictive conditions, $n \geq 42$ only.

{\bf B(7)}. The class $w_{14}$ contains the summand $\langle 8_2 \rangle$. Therefore by (\ref{a1012}) and (\ref{mul}), the product $w_{10} w_{12} w_{14}$ contains the summand $\langle 8_3, 4_5 \rangle $, which is non-trivial if $n \geq 44$.

{\bf C(1).} By (\ref{sw2}), (\ref{b3p}) and (\ref{mul}), 
\begin{equation} \label{c1p} w_2 w_4 w_6 w_8 = \langle 4_2, 2_{14} \rangle + \langle 4_3, 2_{11} \rangle +
\langle 8, 4_2, 2_7 \rangle, \end{equation} which is non-trivial if $n \geq 30$. By Theorem \ref{mthm2}, this proves statement C(1) in the case $N=n+1$, and the case $N=n$ follows by monotonicity.

{\bf C(2)}. By (\ref{sw4}) and (\ref{b5p}), \begin{equation} w_4 w_6 w_8 w_{10} = \langle 8, 4_3, 2_{12} \rangle , \end{equation}
which is non-trivial if $n \geq 44$. This proves statement C(2) for $N=n+3$, which implies it also for $N=n+2$.

{\bf C(3)}. By (\ref{sw6}) and (\ref{b81012})
\begin{equation} \label{c3p} w_6 w_8 w_{10} w_{12} = \langle 8, 4_5, 2_{14} \rangle + \langle 8, 4_7, 2_8 \rangle, \end{equation} which is non-trivial if $n \geq 52$. This proves statement C(3) for $N=n+5$ and hence also for $N=n+4$.

{\bf C(4)}. The class $w_{14}$ contains the summand $\langle 8_2 \rangle$. Therefore by (\ref{b81012}) and (\ref{mul}) the class $w_8 w_{10} w_{12} w_{14}$ contains the summand $\langle 8_3, 4_7, 2_2 \rangle,$ which is non-trivial if $n \geq 56$. This proves statement C(4) for $N=n+7$ and hence also for $N=n+6$.

{\bf D(1)}. By (\ref{c1p}), (\ref{sw10}) and (\ref{mul}), $w_2 w_4 w_6 w_8 w_{10} = \langle 8, 4_3, 2_{14} \rangle,$ which is non-trivial if $n \geq 48$.

{\bf D(2)}. By (\ref{sw4}), (\ref{c3p}) and (\ref{mul}),
$w_4 w_6 w_8 w_{10} w_{12} = \langle 8, 4_7, 2_{12} \rangle , $ which is non-trivial if $n \geq 60$.

{\bf D(3).} Since $w_{14}$ contains summand $\langle 8_2 \rangle$, by (\ref{c3p}) and (\ref{mul}) the product $ w_6 w_8 w_{10} w_{12} w_{14}$ contains summand $\langle 8_3, 4_7, 2_8 \rangle$, which is non-trivial if $n \geq 68$.

{\bf E(1)}. By {\bf D(1)} and formulas (\ref{sw12}) and (\ref{mul}), $w_2 w_4 w_6 w_8 w_{10} w_{12} = \langle 8, 4_7, 2_{14} \rangle$, which is non-trivial if $n \geq 64$.

{\bf E(2)}. Since $w_{14}$ contains the summand $\langle 8_2 \rangle$, by {\bf D(2)} and (\ref{mul}) the product $ w_4 w_6 w_8 w_{10} w_{12} w_{14}$ contains the summand $\langle 8_3, 4_7, 2_{12} \rangle$, which is non-trivial if $n \geq 76$.

{\bf F}. Since $w_{14}$ contains the summand $\langle 8_2 \rangle$, by {\bf E(1)} and formula (\ref{mul}) the class $w_2 w_4 w_6 w_8 w_{10} w_{12} w_{14}$ contains the summand $\langle 8_3, 4_7, 2_{14} \rangle$, non-trivial if $n \geq 80$. \hfill $\diamond$

\section{Equality conditions and homology of knot spaces}

\label{alex}

Denote by $\K$ the affine space of all $C^\infty$-smooth maps $\R^1 \to \R^3$ coinciding with a fixed linear embedding outside some compact set in $\R^1$. Let $\Sigma$ be the {\em discriminant subvariety} of  $\K$  consisting of all maps which are not smooth embeddings, i.e. have either self-intersections or points of vanishing derivative. The elements of the set $\K \setminus \Sigma$ are called {\em long knots}. There is a natural one-to-one correspondence between the connected components of this set and isotopy classes of usual knots, i.e. of smooth embeddings  $S^1 \to \R^3$  or  $S^1 \to S^3$.

The variety $\Sigma$ is swept out by affine subspaces $L(a,b)$ of codimension $3$ in $\K$ corresponding to all possible chords $(a,b)$ in $\R^1$ (including degenerate chords with $a=b$) and consisting of maps  $\varphi: \R^1 \to \R^3$  such that  $\varphi(a) = \varphi(b)$ (or  $\varphi'(a)=0$  if  $a=b$). Much of the topological structure of  $\Sigma$ can be described in terms of the order complex of the (naturally topologized) partially ordered set, whose elements correspond to these subspaces  $L(a,b)$  and their finite intersections (defined by chord diagrams), and the order relation is the incidence. For any $n$, the subspaces in  $\K$  defined in this way by independent $n$-chord diagrams form an affine bundle over the space $\overline{\mbox{CD}}_n$ of equivalence classes of such diagrams (including degenerate ones, containing chords of type $(a,a)$). The fibers of this bundle have codimension $3n$ in $\K$, and its normal bundle is isomorphic to the sum of three copies of the bundle $\tau^*_n$ considered in \S \ref{scheme} (and continued to degenerate chord diagrams).

 The topology of the space $\K \setminus \Sigma$ is related by a kind of Alexander duality to the topology of the complementary space $\Sigma$, in particular, the numerical knot invariants can be realized as linking numbers with infinite-dimensional cycles of codimension 1 in $\K$ contained in $\Sigma$. 
However, Alexander duality deals with finite-dimensional spaces only, therefore to apply it properly we use finite-dimensional approximations of the space $\K$. Namely, we consider infinite sequences  $\K_1 \subset \K_2 \subset \dots $  of finite-dimensional affine subspaces of  $\K$, such that any connected component of  $\K \setminus \Sigma$  is represented by elements of subspaces  $\K_j \setminus \Sigma$  with sufficiently large  $j$, and moreover any homology class of  $\K \setminus \Sigma$  is represented by cycles contained in such subspaces. (The existence of such  sequences of subspaces  $\K_j$  follows easily from Weierstrass approximation theorem). Then for any such subspace  $\K_j$  of dimension  $d_j$ we have Alexander isomorphisms 
\begin{equation}
\label{alexx}
\tilde H^k (\K_j \setminus \Sigma) \simeq \bar H_{d_j-k-1}(\K_j \cap \Sigma),
\end{equation}
 where $\bar H_*$ denotes the Borel--Moore homology groups. 

To study left-hand groups in (\ref{alexx}) (in particular, such a group with  $k=0$, i.e. the group of  $\Z$-valued invariants of knots realizable in  $\K_j$) a {\em simplicial resolution} of the space  $\K_j \cap \Sigma$  is used in \cite{cs}. It is a certain topological space  $\sigma(j)$  and a surjective map $\sigma(j) \to \K_j \cap \Sigma$ inducing an isomorphism of Borel--Moore homology groups. These groups $\bar H_*(\sigma(j)) \simeq \bar H_*(\K_j \cap \Sigma)$ can be calculated by a spectral sequence $\{E^r_{n,\beta}\}$ defined by a natural increasing filtration \begin{equation}
\label{sseq}
\sigma_1(j) \subset \sigma_2(j) \subset \dots \subset \sigma(j), \end{equation}
in particular, $E^1_{n, \beta} \simeq \bar H_{n+\beta}(\sigma_n(j) \setminus \sigma_{n-1}(j))$. 

This filtration is finite if the subspace  $\K_j$  is not very degenerate. Namely, any space $\sigma_n (j) \setminus \sigma_{n-1}(j)$ is  constructed starting from the intersection sets of  $\K_j$  with subspaces of codimension  $3n$  in  $\K$  defined by independent $n$-chord diagrams. Since the family of all such planes is $2n$-parametric, a {\em generic}  $d_j$-dimensional affine subspace $\K_j$ meets only subspaces of this kind with  $n \leq d_j$, so  $\sigma_{d_j}(j)=\sigma(j)$.

\unitlength 1.00mm \linethickness{0.4pt}
\begin{figure}
\begin{picture}(60,75)
\put(0,5){\vector(1,0){80}}
\put(70,5){\vector(0,1){70}}
\put(70,5){\line(-1,1){60}}
\put(78,1){p}
\put(72,73){q}
\put(10,65){\line(0,1){8}}
\put(7,1){$-d_j$}
\put(58,17){\line(0,1){56}}
\put(54.5,0){$\frac{-d_j}{5}$}
\put(20,25){zeros}
\put(30,50){non-stable}
\put(59,40){stable}
\put(68.5,0){$0$}
\end{picture}
\caption{Spectral sequence for $H^*(\K_j \setminus \Sigma)$}
\label{ssf}
\end{figure}
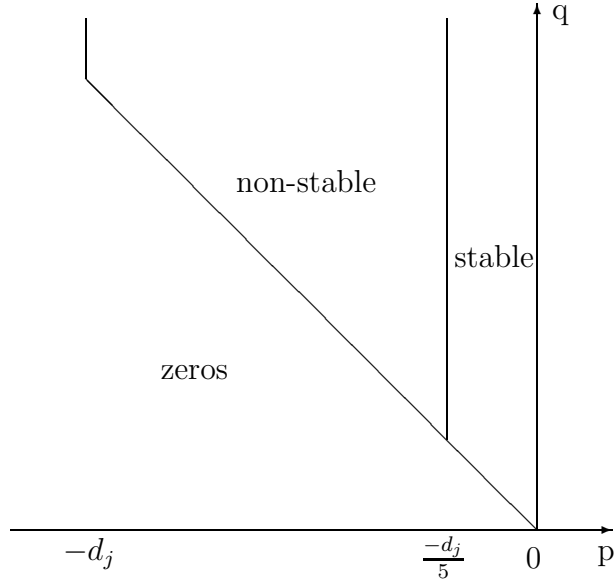

The formal change $E_r^{p,q} \equiv E^r_{-p, d_j-1-q}$ turns the homological spectral sequence defined by this filtration into a cohomological one, which by Alexander duality converges to the left-hand groups of (\ref{alexx}).
 All non-trivial groups $E^{p,q}_r$, $r \geq 1$, of the last spectral sequence for a generic subspace $\K_j$ lie in the domain $\{(p,q): p \in [-d_j,-1], p+q \geq 0\}$, see Fig.~\ref{ssf}.

If the approximating subspace  $\K_j$  is generic and  $n$  is sufficiently small with respect to $d_j$ (namely, $n \leq \frac{d_j}{5}$), then all subspaces of $\K$ defined by independent $n$-chord diagrams intersect $\K_j$ transversally along non-empty planes. Indeed, if  $d_j> 3n $, then the codimension of the set of $d_j$-dimensional affine subspaces in $\K$, which are not generic with respect to a plane of codimension  $3n$,  is equal to  $d_j-3n+1$, therefore the $2n$-parametric family of such sets corresponding to all subspaces defined by $n$-chord diagrams sweeps out a subset of codimension at least $d_j-5n+1$ (if this number is positive), and we can choose $\K_j$ not from this subset.

If  $\K_j$  is generic in this sense, then these intersection sets in  $\K_j$  form an affine bundle of dimension  $d_j -3n$  with base  $\overline{\mbox{CD}}_n$. By the construction of the simplicial resolution, this implies that the topology of sets  $\sigma_n(j) \setminus \sigma_{n-1}(j)$  essentially stabilizes at this value of $j$: 
for any $j' > j$ the space $\sigma_n(j') \setminus \sigma_{n-1}(j')$ is homeomorphic to the direct product of spaces $\sigma_n(j) \setminus \sigma_{n-1}(j)$ and $\R^{d_{j'}-d_j}$. In particular, we have natural isomorphisms $E^1_{n, \beta}(j) \simeq E^1_{n, \beta+(d_j' -d_j)}(j')$ for all $j'>j$, $n \leq d_j/5$  and arbitrary $\beta$. Substitutions (\ref{alexx}) turn them to natural isomorphisms 
$E_1^{p,q}(j') \simeq E_1^{p,q}(j) $ for all  $p \geq - \frac{d_j}{5}$. Moreover, these isomorphisms commute with all further differentials of our spectral sequence, and the Borel--Moore homology groups of the spaces  $\sigma_n(j)$  and  $\sigma_n(j')$   for all  $n \leq \frac{d_j}{5}$  and  $j' \geq j$  are naturally isomorphic to each other up to the shift of dimensions by $d_{j'}-d_j$. The cohomology classes of  $\K \setminus \Sigma$ arising from this area of the spectral sequence (i.e., the sequences of nontrivial cohomology classes of the spaces  $\K_{j'} \setminus \Sigma$, $j'\geq j$, realizable by linking numbers with cycles lying in $\sigma_n(j')$, $n \leq \frac{d_j}{5}$, and corresponding to one another by these isomorphisms) are known as {\em finite type} cohomology classes of the space of knots. Therefore, the intriguing question on the completeness of the system of these classes in entire cohomology groups of  $\K \setminus \Sigma$
(in particular, on the existence of non-equivalent knots not separated by finite type invariants) depends on the groups $E_r^{p,q}(j)$ in the non-stable domains, on the deviation of these groups from stable ones, and on the way in which the non-stable groups $E^{p,q}_\infty(j)$ for different  $j$   correspond to the same cohomology classes of spaces  $\K_j \setminus \Sigma$  with different  $j$. 

The arguments of the previous sections of this article allow us to say something about the non-triviality of this problem.

\begin{proposition}
\label{proapp1}
If  $4n-I(n) > d_j \geq 3n$,  then for any  $d_j$-dimensional affine subspace  $\K_j \subset \K$  there exist independent $n$-chord diagrams such that corresponding affine subspaces of  $\K$  have non-generic $($i.e. either non-transversal or empty$)$
intersections with the space  $\K_j$.
\end{proposition}

\begin{proposition}
\label{proapp2}
If  $2n+I(n) \leq d_j \leq 3n$, then for almost any  $d_j$-dimensional affine subspace  $\K_j \subset \K$  $($i.e., for any subspace from a residual subset in the space of all such subspaces$)$ there exist independent  $n$-chord diagrams such that the corresponding affine subspaces of  $\K$  have non-empty intersection with  $\K_j$.
\end{proposition}

\begin{definition} \rm 
Denote by  ${\LL}$  the affine bundle over the space  $B(\R^2_+,n) \setminus \Xi$ \
of independent  $n$-chord diagrams, whose fiber over any such diagram is the subspace of codimension  $3n$  in  $\K$ consisting of maps  $\varphi: \R^1 \to \R^3$  taking the same values at endpoints of each chord of this diagram.
For an affine subspace  $\K_j \subset \K$  denote by  $\|(\K_j)$  the subset in  $B(\R^2_+, n) \setminus \Xi$  consisting of  $n$-chord diagrams such that the corresponding fiber of bundle  $\LL$  contains lines parallel to some lines lying in the space  $\K_j$.
\end{definition}

\noindent 
{\it Proof of Proposition \ref{proapp1}.}  The normal bundle  ${\LL}^\perp$  of   $\LL$  in  $\K$  is isomorphic to the direct sum of three copies of the regular bundle  $\xi_n$. By Lemmas \ref{le1} and \ref{le2}, its total Stiefel--Whitney class is then equal to $(w(\xi_n))^3 \equiv w(\xi_n),$ in particular, its $i$-dimensional component  $w_i$  is not trivial if  $i \leq n-I(n)$. 

Make  $\K_j$  a vector space by choosing arbitrarily the ``origin'' point in it. If all fibers of the bundle  $\LL$  are in general position with respect to  $\K_j$, then a $(d_j-3n)$-dimensional vector bundle with the same base is defined, the fiber of which over a chord diagram is obtained from the intersection set of  $\K_j$  and the corresponding fiber of the bundle $\LL$ by a parallel translation, after which it passes through the origin point of $\K_j$. The total Stiefel--Whitney class of this bundle is equal to  $w({\LL}^\perp)^{-1} \equiv w(\xi_n)^{-1}$, which by Lemma \ref{le2} is equal to  $w(\xi_n)$. If  $d_j -3n < n-I(n)$  then this implies that  $w_{n-I(n)}(\xi_n)=0$, a contradiction. \hfill $\diamond$

\begin{lemma}
If  $d_j \leq 3n$, then for almost any  $d_j$-dimensional affine subspace  $\K_j \subset \K$ the codimension of the set  $\|(\K_j)$  in  $B(\R^2_+, n) \setminus \Xi$  is at least  $3n-d_j+1$. 
\end{lemma}

\noindent
{\it Proof.} Consider the space 
\begin{equation}
\label{tb}
\tilde G(\K, d_j) \times (B(\R^2_+,2) \setminus \Xi)
\end{equation}
of all pairs $\{\K_j, \Gamma\}$ where $K_j$ is a $d_j$-dimensional affine subspace of  $\K$ and $\Gamma$ is an independent $n$-chord diagram. Denote by $\Lambda$ the subset of this space consisting of pairs  $\{\K_j, \Gamma\}$  such   that $\Gamma \in \|(K_j)$.
Both the space (\ref{tb}) and its subset $\Lambda$ are fibered over the space $B(\R^2_+,2) \setminus \Xi$  of independent  $n$-chord diagrams, and for any such diagram $\Gamma$ the corresponding fiber of the latter fiber bundle has codimension $3n-d_j+1$ in the fiber of the former. Therefore, the codimension of $\Lambda$ in the space (\ref{tb}) is equal to  $3n-d_j+1$, and the typical fiber of the projection of $\Lambda$ to the first factor of (\ref{tb})  has codimension at least  $3n-d_j+1$ in the corresponding fiber  of the projection of entire space (\ref{tb}). \hfill $\diamond$ \medskip

\noindent
{\it Proof of Proposition \ref{proapp2}.}
Let us fix a subspace  $\K_j$  for which the condition of the previous lemma is satisfied. The complement of the set  $\|(\K_j)$  in the manifold  $B(\R^2_+,n) \setminus \Xi$  has then the same homology groups up to dimension  $3n-d_j$  as entire  $B(\R^2_+, n) \setminus \Xi$. 

Consider the affine bundle  $(\LL^\perp)^*$  over the manifold $B(\R^2_+,n) \setminus \Xi$ : its fibers consist of linear functions on  $\K$  vanishing on the corresponding fibers of the bundle $\LL$. Over the set $(B(\R^2_+,n) \setminus \Xi) \setminus \|(K_j)$ a $(3n-d_j)$-dimensional subbundle of $(\LL^\perp)^*$ is defined, consisting of functions constant on  $\K_j$. This subbundle has the same Stiefel--Whitney class (equal to $w(\xi_n)$) as the whole bundle $(\LL^\perp)^*$, since its normal bundle is isomorphic to the trivial bundle with fiber $(\K_j)^*$. If no fibers of the bundle  $\LL$  intersect the space  $\K_j$, then this subbundle has a nowhere vanishing cross-section: indeed, we can define an arbitrary Euclidean structure on this subbundle, and choose in each fiber the linear function of unit norm taking the maximal value on  $\K_j$. If  $3n-d_j \leq n-I(n)$  then this contradicts the non-triviality of the class  $w_{n-I(n)}(\xi_n)$. \hfill $\diamond$

\begin{remark} \rm
I hope that the further study of characteristic classes of the bundle  $\LL$  (and of its analog defined on the space  $\overline{CD}_n$  of equivalence classes of chord diagrams, rather than on the resolution  $B(\R^2_+,n) \setminus \Xi$  of this space) will provide not only the proofs of unavoidable troubles in calculating the cohomology classes of knot spaces, but also the construction of some such classes not reducible to classes of finite type.
\end{remark}

\subsection*{Acknowledgments}

I thank very much Weizmann Institute, where this work was completed, and especially Sergey Yakovenko for extremely timely help and hospitality.
I thank M.E.~Kazarian for a consultation on the Thom--Porteous formula.


\begin{thebibliography}{9}

\bibitem{CH} F.R.~ Cohen and D.~Handel, {\it $k$-regular embeddings of the plane},
Proceedings of the American Mathematical Society 72 (1) (1978), 201-204.

\bibitem{CC} {\it Configuration spaces}. Geometry, topology and representation theory. Springer, 2016. xii+379 pp.

\bibitem{fuks} D.B.~Fuchs, {\it Cohomologies of the braid group mod 2}, Funct. Anal. Appl., 4:2 (1970), 143–-151.

\bibitem{ful} W.~Fulton, {\it Intersection Theory}, Springer, 1984.

%\bibitem{GM} M.~Goresky, R.~MacPherson, {\it Stratified Morse Theory}, Springer, Berlin-Heidelberg, 1988.

\bibitem{MS} J.W.~Milnor and J.D.~Stasheff, {\it Characteristic classes}, Princeton University Press, Princeton, N. J.; University of Tokyo Press, Tokyo, 1974. Annals of Mathematics Studies, No. 76.

%\bibitem{port} I.R.~Porteous, {\it Simple singularities of maps}. In: Proceedings of Liverpool Singularities--Symposium I. Lect. Notes in Math., vol. 192 (1971), 286--307, Springer, Berlin--Heidelberg.

\bibitem{cs} V.A.~Vassiliev, {\it Cohomology of knot spaces}. In: Theory of singularities and its applications, Adv. Soviet Math., 1, American Mathematical Society, Providence, RI, 1990, 23–-69.

\end{thebibliography}
\end{document}